\documentclass[12pt]{amsart}
\usepackage[all]{xy}
\xyoption{import}
\usepackage{graphicx,amssymb,rotating,hyperref}
\advance\hoffset by -0.5 truecm
\newtheorem{Theorem}{Theorem}[section]
\newtheorem{Lemma}[Theorem]{Lemma}

\newtheorem{Proposition}[Theorem]{Proposition}

\newtheorem{Remark}[Theorem]{Remark}

\def\ds{{\bf ds}}

\def\Z{{\mathbb Z}}
\def\R{{\mathbb R}}

\def\f{\tilde{f}_{k}}

\def \re{{\mathbb R}}
\def \Q{{\mathbb Q}}

\def \0{\lambda_{0}}

\def \S{\mbox{\bf Sol}}
\def \g{{\mathfrak g}}

\def\u{{u}}
\def\y{{y}}
\def\a{{\alpha}}
\def\n{{\nu}}
\def\s{{\it s}}
\def\Sol{{\bf Sol}}
\def\sol{{\mathfrak s}}
\def\d{{\rm d}}
\def\G{{\rm G}}
\def\Ss{{\bf S}}
\def\bn{{\bar\n}}
\def\f{{\rm f}}
\def\m{{\tt m}}
\def\leb{{{\pmb{\mu}}}}
\def\DSol{{\pmb{\Delta}}}
\newcommand\pb[2] {{\left\{#1,#2\right\}_{\s}}}
\newcommand\vol[1] {{{\rm vol}\, #1 }}
\def\T{{\bf T}}
\def\O{{\mathfrak O}}
\def\U{{\mathfrak U}}
\def\cc{\gamma}

\begin{document}
\title[Magnetic flows on $\S$-manifolds]{Magnetic flows on $\S$-manifolds:
dynamical and symplectic aspects}

\author[L.T. Butler]{Leo T. Butler}
\address{School of Mathematics, The University of Edinburgh, 6214 James Clerk
Maxwell building, Edinburgh, EH9 3JZ, UK}
\email{l.butler@ed.ac.uk}

\author[G.P. Paternain]{Gabriel P. Paternain}
 \address{ Department of Pure Mathematics and Mathematical Statistics,
University of Cambridge,
Cambridge CB3 0WB, UK}
 \email {g.p.paternain@dpmms.cam.ac.uk}




\begin{abstract} We consider magnetic flows on compact quotients
of the 3-dimensional solvable geometry $\S$ determined by the
usual left-invariant metric and the distinguished monopole. We show that
these flows have positive Liouville entropy and therefore
are never completely integrable. 
This should be compared with the known fact that the underlying
geodesic flow is completely integrable in spite of having positive
topological entropy.
We also show that for a large class of twisted cotangent bundles of 
solvable manifolds every compact set is displaceable.

\end{abstract}

\maketitle

\section{Introduction}The Lie group $\S$ is the semidirect product
associated with the
action of $\re$ on $\re^2$ given by
\[u \cdot \,(y_0,y_1)=(e^u y_0, e^{-u} y_1).\] 
The group $\S$ is diffeomorphic to $\re^3$ and the product is
\[(y_0,y_1,u)\star(y_{0}',y_{1}',u')=(e^u y_{0}'+y_{0},e^{-u} y_{1}'+y_{1},u+u').\]
It is not difficult to see that $\S$ admits cocompact lattices.
Let $A\in SL(2,\mathbb Z)$ be such that there is $P\in GL(2,\mathbb R)$ with
\[PAP^{-1}=\left(\begin{array}{cc}
\lambda &0\\
0&1/\lambda\\
\end{array}\right)\]
and $\lambda>1$.
There is an injective homomorphism
$$\mathbb Z^2\ltimes_{A}\mathbb{Z}\hookrightarrow \mathbf{Sol}$$
given by $(m,n,l)\mapsto (P(m,n),\log\lambda\,l)$ which defines
a cocompact lattice $\DSol$ in $\mathbf{Sol}$.
The closed 3-manifold $\Sigma:=\DSol\setminus \S$ is a 2-torus
bundle over the circle with hyperbolic gluing map $A$.

The Riemannian metric
\[\ds^2=e^{-2u}dy_{0}^2+e^{2u}dy_{1}^2+du^2\]
is left-invariant and descends to a Riemannian metric on $\Sigma$.
It is a remarkable fact discovered by A. Bolsinov and I. Taimanov \cite{BT}
that the geodesic flow of $(\Sigma,\ds^2)$ is completely integrable in the sense
of Liouville with the two additional integrals
\begin{align*}
f&=p_{y_{0}}p_{y_{1}}\\
F&=\exp\left(\frac{-1}{p^{2}_{y_{0}}p_{y_{1}}^2}\right)\sin\left(2\pi\frac{\log|p_{y_{0}}|}{\log\lambda}\right).\\
\end{align*}
The geodesic flow has topological entropy $h_{top}=1$ but Liouville
(or metric) entropy $h_{\leb}=0$. It is the simplest example of a
geodesic flow on
a compact homogeneous space with these properties.  Note that the
lattice $\DSol$ has exponential word growth and the entropy is all
carried in the minimizing Aubry-Mather sets given by $p_{u}=\pm 1,
p_{y_0}=p_{y_1}=0$.  The dynamics on these sets is Anosov and given by
the suspension of $A$. We refer to \cite{BDV} for a detailed
description of the foliation by Liouville tori and for spectral properties
of the Laplace-Beltrami operator of $(\Sigma,\ds^2)$.

The manifold $\Sigma$ has a distinguished {\it monopole}, i.e. a closed
non-exact 2-form which generates $H^{2}(\Sigma,\mathbb R)$ given by
$\Omega=dy_{0}\wedge dy_{1}$. This form is harmonic and Hodge dual to
the generator $du$ of $H^{1}(\Sigma,\mathbb R)$. The Aubry-Mather sets we
mentioned before are calibrated by the closed 1-forms $\pm du$.

The first goal of this paper is the study of the dynamics of the
magnetic flow determined by the metric $\ds^2$ and the monopole
$\Omega$.  We will modulate the intensity of the magnetic field
$\Omega$ with a parameter $s\in [0,\infty)$ and we will always
  consider the magnetic flow $\varphi^s$ running with speed one. The
  analysis of the flow is carried out in Section \ref{amf}. One of our
  findings is that the magnetic flow ceases to be Liouville integrable
  as soon as the magnetic field is switched on.  The reason is that
  the Liouville entropy becomes positive. In fact, one can compute the
  Liouville entropy exactly as we now explain.  Since all the objects
  involved are left-invariant the flow $\varphi^s$ may be reduced to
  an {\it Euler} flow $\psi^s$ on $\mathfrak s^*$, the dual of the Lie
  algebra $\mathfrak s$ of $\S$.  With respect to the basis of
  left-invariant 1-forms $\{e^{-u}dy_{0}, e^{u}dy_{1}, du\}$, a point
  in $\mathfrak s^{*}$ will have coordinates $(\a_{0},\a_{1},\n)$.  It
  is easy to see that $f=\a_{0}\a_{1}+s\n$ is a Casimir and thus an
  integral of $\psi^s$. Observe that $\psi^s$ leaves invariant the
  sphere $\Ss$ given by $\a_{0}^2+\a_{1}^2+\n^2=1$. Let $d\theta$ be
  its canonical probability area measure.

\medskip

\noindent {\bf Theorem A.} {\it The Liouville entropy of $\varphi^s$
is given by\[h_{\leb}(\varphi^{s})=\int_{\Ss}|\bar{\n}|\,d\theta\]
where $\bar{\n}$ is the average of $\n$ over the level sets of the
Casimir $f$. Moreover, $h_{\leb}(\varphi^{s})>0$ for all $s>0$ and
approaches $1/2$ as $s\to\infty$, while $h_{top}(\varphi^s)\equiv 1$.
}

\medskip

This result should be compared with the well-known example of the
magnetic flow on a compact hyperbolic surface with magnetic field
given by the area form. In this example, as the intensity $s$
increases the flow becomes ``simpler''. Indeed, topological entropy
decreases; at $s=1$ we hit the horocycle flow and for $s>1$, the flow
has all its orbits closed and becomes integrable.  The opposite seems
to be happening for our magnetic flow on $\S$. On the other hand, the
well-known Rydberg model of a hydrogen atom in a strong magnetic field
is believed to exhibit behaviour similar to that described in Theorem
A. We are unaware of any {\em proof}, as opposed to evidence, that the
Rydberg model has positive Liouville entropy.

The second goal of this paper is to try to explain these drastic changes
in the dynamics in terms of changes in the symplectic topology of twisted
cotangent bundles.

Let $\Sigma$ be a closed manifold and let $\omega_0$ be the canonical
symplectic form of the cotangent bundle $\tau:T^*\Sigma\to \Sigma$.
Given a closed 2-form $\sigma$ we let
$\omega_{\sigma}:=\omega_{0}-\tau^{*}\sigma$ be the twisted symplectic
form determined by $\sigma$.  Recall that given a compact set $K$, the
{\it displacement energy} of $K$ is defined as
\[e(K):=\inf \{\rho(\mathbf{1},h):\;h\in \mbox{\rm Ham}_{c}(T^{*}\Sigma,\omega_{\sigma}),\;\;
h(K)\cap K=\emptyset\}\] where $\rho$ is Hofer's distance and
$\mbox{\rm Ham}_{c}(T^{*}\Sigma,\omega_{\sigma})$ is the set of
compactly supported Hamiltonian diffeomorphisms.  Recall also that a
compact set $K$ is said to be {\it displaceable} if there exists $h\in
\mbox{\rm Ham}_{c}(T^{*}\Sigma,\omega_{\sigma})$ such that $h(K)\cap
K=\emptyset$. Thus $K$ is displaceable iff $e(K)$ is finite. A well
known result of M. Gromov \cite{G} asserts that the zero section of
$(T^{*}\Sigma,\omega_{0})$ is {\it not} displaceable. On the other
hand, if $\sigma$ is non-zero and $\Sigma$ has zero Euler
characteristic, results of F. Laudenbach and J.-C. Sikorav \cite{LS}
and L. Polterovich \cite{Po} imply that the zero section of
$(T^{*}\Sigma,\omega_{\sigma})$ is actually displaceable (if $\sigma$
is non-zero, the zero section of $T^*\Sigma$ ceases to be Lagrangian).
Finite displacement energy has important implications. According to a
recent result of F. Schlenk \cite{S}, if a compact energy level of an
autonomous Hamiltonian is displaceable, then it will have finite
$\pi_1$-sensitive Hofer-Zehnder capacity which in turn yields almost
everywhere existence of {\it contractible} closed orbits (i.e. there is
a full measure set of values of the energy for which the corresponding
energy level has a contractible closed orbit). Let us
illustrate this discussion with the following example. Consider a
closed hyperbolic 3-manifold and let $\sigma$ be {\it any} non-zero
closed 2-form. For high values of the energy the magnetic flow will be
Anosov, since it can be seen as a pertubation of a geodesic flow on a
negatively curved manifold.  Thus for high energies, the magnetic flow
will have no contractible closed orbits (the magnetic flow will be
topologically conjugate to the geodesic flow and it is well known that
the latter has no contractible closed geodesics). Schlenk's result now
implies that high energy levels are not displaceable, while low energy
levels are by the results of Laudenbach-Sikorav and Polterovich. If we
take the closed 3-manifold to have non-zero first Betti number, then it will
have non-zero second Betti number and
we may choose magnetic fields $\sigma$ with non-zero cohomology
classes (monopoles).

Returning to our example on $\S$ we note that the geodesic flow of
$(\Sigma, \ds^2)$ has no contractible closed orbits, but as soon as the
magnetic field is switched on, contractible closed orbits
appear. These orbits are related to the vanishing of $\bar{\n}$, see
Remark \ref{clor} were these observations are proved. 
It turns out that {\it every} compact set in
$(T^*(\DSol\setminus\S), \omega_{\Omega})$ is displaceable. Our
last result shows that this is also true for a large class of solvable
manifolds.

We say that a Lie group $G$ is {\it completely solvable} if it is a
closed subgroup of the group of upper triangular matrices with
positive diagonal entries. The class of completely solvable groups
lies strictly in between nilpotent and solvable groups.

Given a Lie algebra $\mathfrak g$, let $L:\Lambda_{2}(\mathfrak
g)\to\mathfrak g$ be the linear map induced by the Lie bracket, where
$\Lambda_{2}(\mathfrak g)$ is the second exterior power of $\mathfrak g$.
 Recall
that 2-vectors are elements in $\Lambda_{2}(\mathfrak g)$ of the form
$x\wedge y$ with $x,y\in \mathfrak g$.

\medskip

\noindent {\bf Theorem B.} {\it Let $G$ be a simply connected
completely solvable group and suppose $\mbox{\rm Ker}\,L$ is generated
by 2-vectors. Let $\Gamma$ be a cocompact lattice and
$\Sigma:=\Gamma\setminus G$.  Then, for any monopole $\sigma$ and any
compact set $K\subset (T^{*}\Sigma,\omega_{\sigma})$, $e(K)<\infty$.
}

\medskip

Certainly, our example $(T^*(\DSol\setminus\S),\omega_{\Omega})$
fits the hypotheses of the theorem. For tori,
the theorem also follows from the proof of Theorem 3.1 in \cite{GK}. 
 It is quite likely that Theorem B
holds for any simply connected solvable Lie group with lattice.  We do
not know of an example of a solvable Lie algebra where $\mbox{\rm
Ker}\,L$ is not generated by 2-vectors. In Section \ref{tb} we show
how Theorem B applies to compact quotients of some of the standard
nilpotent Lie algebras, like the Heisenberg Lie algebra $\mathfrak
h_{2n+1}$ and the Lie algebra of upper triangular matrices $\mathfrak
u_{n}$.  Finally, in Subsection \ref{mcv} we discuss these results in
the context of Aubry-Mather theory and Ma\~{n}\'e's critical values.

\medskip

{\it Acknowledgement:} We would like to thank L. Polterovich and F. Schlenk
for useful comments and discussions about Theorem B. 

\section{Preliminaries}
Let $\Sol$ be the semidirect product of $\R^2$ with $\R$, with
coordinates $(\u,\y_0,\y_1)$ and multiplication
\begin{equation} \label{eq:sm}
(\y_0,\y_1,u)\star(\y_0',\y_1',u')=(\y_0+e^{\u}\y_0', \y_1+e^{-\u}\y_1',u+u').
\end{equation}
The map $(\y_0,\y_1,u) \mapsto \u$ is the epimorphism $\Sol \to \R$
whose kernel is the normal subgroup $\R^2$.
The group $\S$ is isomorphic to the matrix group
\[\left(\begin{array}{ccc}

e^u&0&y_0\\
0&e^{-u}&y_1\\
0&0&1\\

\end{array}\right).\]

If one denotes by $p_{\u}$, $p_{\y_0}$ and $p_{\y_1}$ the momenta that
are canonically conjugate to $\u$, $\y_0$ and $\y_1$ respectively,
then the functions
\begin{equation} \label{eq:mom}
\begin{array}{lcl}
\a_0 &=& e^{\u} p_{\y_0}, \\
\a_1 &=& e^{-\u} p_{\y_1}, \\
\n &=& p_{\u}
\end{array}
\end{equation}
are left-invariant functions on $T^*\Sol$. The closed $2$-form
\begin{equation} \label{eq:magnetic}
\Omega = \d\y_0 \wedge \d\y_1
\end{equation}
is also left-invariant, and consequently,
\begin{equation} \label{eq:sympform}
\omega_{\s} = \d p_{\u} \wedge \d\u + \d p_{\y_0} \wedge \d\y_0 +
\d p_{\y_1} \wedge \d\y_1 - \s\d\y_0 \wedge \d\y_1 
\end{equation}
is a left-invariant twisted symplectic form on $T^*\Sol$ for any real number
$\s$. The Poisson bracket induced by $\omega_{\s}$ is denoted by
$\pb{}{}$. The Poisson brackets of the coordinate functions are
\begin{equation} \label{eq:pb}
\begin{array}{lclclcl}
\pb{\n}{\u} &=& 1,      &\hspace{10mm}& \pb{\a_0}{\a_1} &=& s,\\
\pb{\a_0}{\y_0} &=& e^{\u},  && \pb{\n}{\a_0} &=& \a_0, \\
\pb{\a_1}{\y_1} &=& e^{-\u},  && \pb{\n}{\a_1} &=& -\a_1,
\end{array}
\end{equation}
and all others vanish. Define the Hamiltonian $H$ on $T^*\Sol$ by
\begin{equation} \label{eq:H}
2H = \n^2 + \a_0^2 + \a_1^2,
\end{equation}
so that when $\s=0$, $H$ is the Hamiltonian of the left-invariant
Riemannian metric mentioned in the Introduction. The equations of the magnetic flow induced by
$H$ are
\begin{equation} \label{eq:XH}
X_H = \left\{ \begin{array}{lclclcl}
\dot{\u}   &=& \n,          & \hspace{10mm} & \dot{\n}   &=& -\a_0^2 + \a_1^2,\\
\dot{\y}_0 &=& e^{\u} \a_0, &            & \dot{\a}_0 &=& -\a_1 s + \n\a_0,\\
\dot{\y}_1 &=& e^{-\u} \a_1, &           & \dot{\a}_1 &=&  \a_0 s - \n\a_1,
\end{array}
\right.
\end{equation}
or $X_H(\bullet) = \pb{H}{ \bullet }$.

The Lie algebra of left-invariant functions on $T^*\Sol$ has a
non-trivial centre generated by the Casimir
\begin{equation} \label{eq:cas}
f = s\n + \a_0 \a_1.
\end{equation}

\begin{Remark}{\rm The 2-form $\Omega$ defines a central extension
of $\Sol$: $\R \hookrightarrow \G \to \Sol$. The Lie algebra $\g$ of
$\G$ is isomorphic to the Lie algebra with basis $s,\n,\a_0,\a_1$ and
Lie bracket $\pb{}{}$. The equations of the magnetic Hamiltonian $H$
(equation \ref{eq:XH}) may be viewed as the symplectic reduction of a
Kaluza-Klein metric Hamiltonian on $T^*\G$ at a non-zero level of
momentum. From this point of view, $f$ and $s$ are Casimirs of the
Poisson bracket on $\g^*$.

Actually, the group $\G$ may be identified with one of the solvable
4-dimensional geometries, namely $\mbox{\rm Sol}^{4}_{1}$ \cite{W}.
It has a matrix representation
\[\left(\begin{array}{ccc}

1&y&z\\
0&e^{t}&x\\
0&0&1\\

\end{array}\right),\]
where $x,y,z,t\in\R$. Via the Kaluza-Klein metric, Theorem A could be 
reinterpreted as follows: the geodesic flow on compact quotients
of $\mbox{\rm Sol}^{4}_{1}$ has positive Liouville entropy and is not
completely integrable.
 }
\end{Remark}

\section{Analysis of the Magnetic Flow}
\label{amf}
Since the Hamiltonian vector field $X_H$ (equation \ref{eq:XH}) is
left-invariant, the vector field factors onto a vector field $E_h$ on
$\sol^*$ through the projection map $T^*\Sol \to \sol^*$ induced by
the left-framing of $T^*\Sol$. The {\em Euler} vector field $E_h$ is a
Hamiltonian vector field on $\sol^*$ equipped with the Lie bracket
$\pb{}{}$. The Hamiltonian $h : \sol^* \to \R$ is the Hamiltonian
which induces $H$. It is clear the dynamics of $X_H$ can be
reconstructed from the dynamics of $E_h$.

Let $\Ss = h^{-1}(\frac{1}{2})$ be the unit sphere in $\sol^*$; the
unit-sphere bundle $H^{-1}(\frac{1}{2})$ is naturally diffeomorphic to
$\Sol \times \Ss$. The functions $\n,\a_0,\a_1$ will be regarded as
coordinate functions on $\sol^*$. Define the standard smooth measure
$\theta$ on $\Ss$ by
\begin{equation} \label{eq:mtheta}
4\pi \times \theta = \left. \n \d\a_0\wedge \d\a_1 + \a_0 \d\a_1\wedge \d\n +
\a_1 \d\n \wedge \d\a_0 \right|_{\Ss}.
\end{equation}
The measure $\theta$ may be decomposed as $\theta = \m \wedge
\m_f$. The measure $\m$ is defined so that for each connected
component of $f^{-1}(c) \cap \Ss$, call it $\f_c$, $\m$ induces a
smooth probability measure on $\f_c$ that is $E_h$-invariant. Let $\bn
: \Ss \to \R$ be defined by
\begin{equation} \label{eq:bn}
\bn(\mu) := \oint_{\f_{f(\mu)}}\ \n\, \d \m \hspace{20mm} \forall
\mu \in \Ss,
\end{equation}
that is, $\bn(\mu)$ is the mean value of $\n$ along the connected
component of the level set of $f|\Ss$ containing $\mu$.

Here is a more prosaic definition of $\m$. Because the vector field
$E_h$ preserves the volume form $\d \n \wedge \d \a_0 \wedge \d \a_1$
on $\sol^*$, and $E_h$ is tangent to the unit sphere $\Ss$, the vector
field $E_h | \Ss$ is Hamiltonian with respect to the symplectic form
$\theta$ (the Hamiltonian is $g = 4\pi \times f$). Therefore, if $c$ is a
non-trivial regular value of the integral $f$, then a neighbourhood of
$\f_c$ in $\Ss$ admits action-angle coordinates $(I,\phi \bmod 1)$
such that $g=g(I)$,
\begin{equation} \label{eq:EhI}
E_h = \left\{
\begin{array}{lcl}
\dot{\phi} &=& \frac{\partial g(I) }{\partial I},\\
\dot{I} &=& 0,
\end{array}
\right.
\end{equation}
and $\theta = \d\phi \wedge \d I$. The measure $\m$ in these
coordinates is
\begin{equation} \label{eq:m}
\m = \d \phi,
\end{equation}
while
\begin{equation} \label{eq:bnloc}
\bn = \int_0^1 \n(\phi,I)\ \d \phi.
\end{equation}

\begin{Proposition}
For $s\neq 0$, $\bn : \Ss\to\R$
  is a continuous, $\psi^s$-invariant function which is real-analytic
  off the set of non-elliptic singular levels of $f|\Ss$.
\label{PA}
\end{Proposition}

\begin{center}
\begin{figure}[!ht]
\includegraphics[width=8cm, height=8cm]{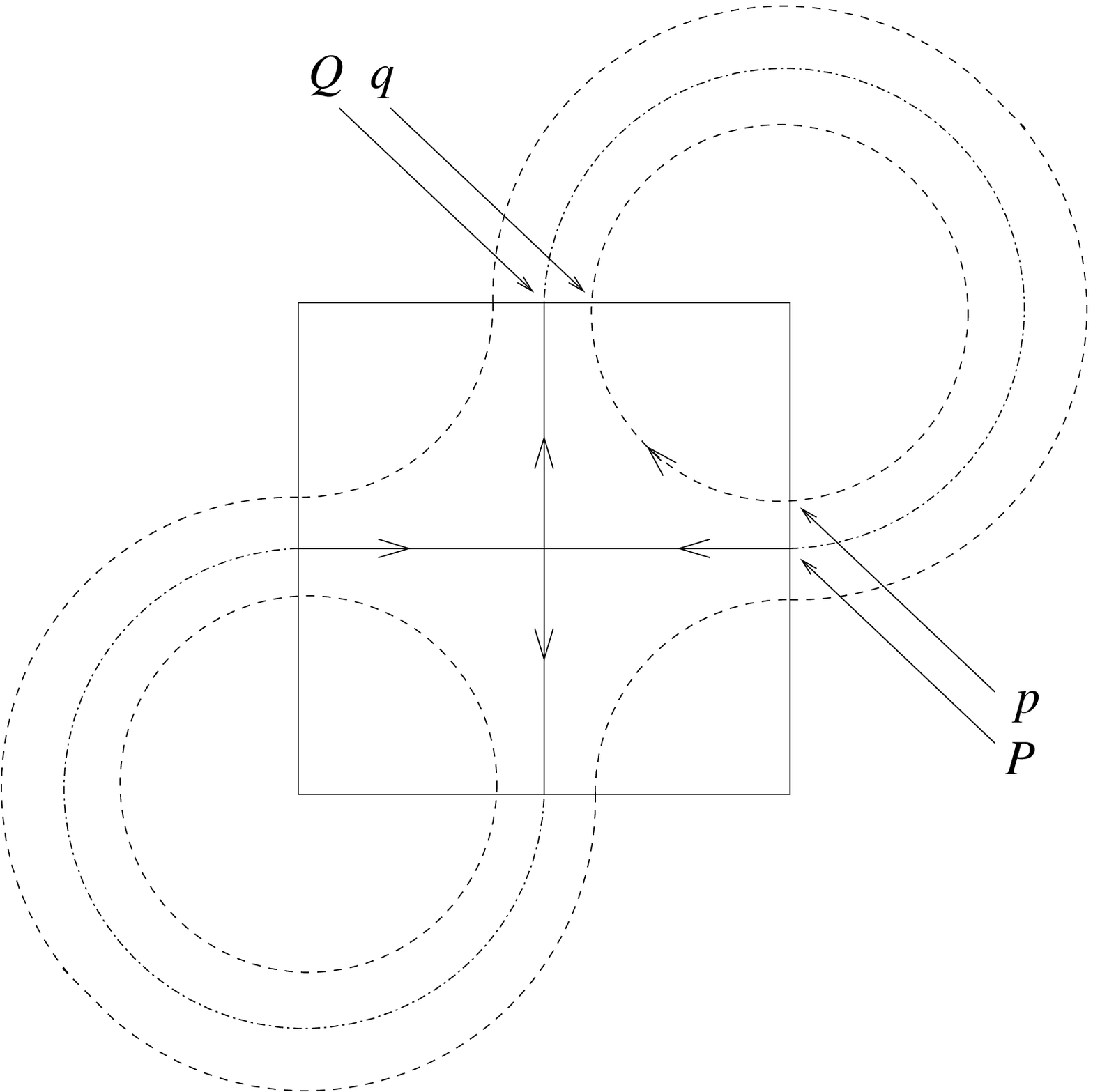}
\caption{} \label{fig:2}
\end{figure}
\end{center}

\begin{proof}
The real-analyticity of $\bn$ on the regular-point set follows from
the fact that $\bn$ and $f$ are real-analytic and the action-angle
coordinates are real-analytic.

\noindent{\em Case 1, $|s| \neq 0,1$:} When $|s|<1$, $f$ has a pair of
peaks (resp. pits) at $\a_0=\a_1=\pm\a,\n=s$
(resp. $\a_0=-\a_1=\pm\a,\n=-s$) where
$\a=\sqrt{\frac{1}{2}(1-s^2)}$. When $|s|\geq 1$, $f$ has a single
peak (resp. pit) at $\a_0=\a_1=0,\n=1$
(resp. $\a_0=\a_1=0,\n=-1$). These critical points are all
non-degenerate for $|s| \neq 0,1$.

\begin{center}
\begin{figure}[htb]
\includegraphics[width=5cm, height=4cm]{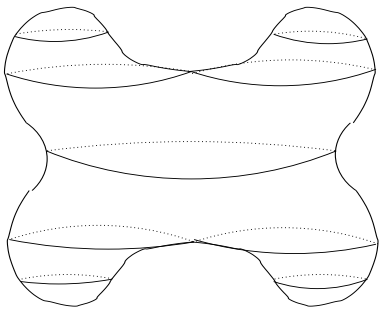}
\caption{$\Ss$ seen from the point of view of $f$, $0<|s|<1$.} \label{fig:1}
\end{figure}
\end{center}

\medskip
\noindent{\em Case 1a, elliptic singularity:} Let $p \in \Ss$ be a
peak or pit for $f|\Ss$, hence an elliptic singularity of $E_h$ on
$\Ss$. There is a canonical system of coordinates $(x,y)$ defined on a
neighbourhood of $p$ such that the Hamiltonian $g$ of $E_h|\Ss$ is in
Birkhoff normal form:
\begin{equation} \label{eq:gbnf}
g(x,y) = g_1I + g_2 I^2 + \cdots, \hspace{10mm} I=\frac{1}{2}\left( x^2
+ y^2 \right), x+iy = \sqrt{2I} e^{2\pi i \phi}.
\end{equation}
It is well-known that $g$ has a formal Birkhoff normal form; Zung has
proven that the formal Birkhoff normal form converges when $g$ is
completely integrable~\cite{Z}. Inspection of
equations~(\ref{eq:EhI}--\ref{eq:bnloc}) shows that $\bn$ may be
written as
\begin{equation} \label{eq:bnloc2}
\bn(\mu) = \frac{1}{T} \times \int_0^T \n \circ \psi^s_t(\mu)\, \d t,
\hspace{10mm} \forall \mu \in \Ss,
\end{equation}
where $\psi^s$ is the Euler flow of $E_h|\Ss$ and $T$ is the period of
the orbit through $\mu$. In an action-angle chart $T=\frac{\partial
I}{\partial g}$, and one sees that $T$ extends over the critical point
at $I=0$ as a real-analytic function. Therefore, define
\begin{equation} \label{eq:s1}
t\cdot \mu = \psi^s_{tT(\mu)}(\mu), \hspace{10mm} \forall t \in S^1=\R/\Z.
\end{equation}
This defines a real-analytic action of $S^1$ on a neighbourhood of the
critical point $p$. In angle-action coordinates, this action is just
$t \cdot (\phi,I) = (\phi+t \bmod 1, I)$. The integral in equation
(\ref{eq:bnloc2}) is then
\begin{equation} \label{eq:bnloc2-1}
\bn(\mu) = \int_0^1 \n(t \cdot \mu)\, \d t, \hspace{10mm} \forall \mu\in\Ss
\end{equation}
{\it i.e.} $\bn$ is the average of $\n$ under the real-analytic action
of $S^1$. This shows that $\bn$ is real-analytic in a neighbourhood of
the elliptic critical point $p$.

\medskip
\noindent{\em Case 1b, hyperbolic singularity:} In this case,
it is known that there are canonical coordinates $(x,y)$
which send the hyperbolic fixed point to $(0,0)$, its stable and
unstable manifolds to the $x$- and $y$-axes respectively, and in which
the hamiltonian is of the form
\begin{equation} \label{eq:ghyp}
g = g_1 \tau + g_2 \tau^2 + \cdots, \hspace{10mm} \textrm{where}\ \tau=xy.
\end{equation}
In this coordinate system, the flow is simply
\begin{equation} \label{eq:psihyp}
\psi^s_t(x,y) = (x e^{-t \omega(\tau)}, y^{t \omega(\tau)})
\end{equation}
where $\omega =\frac{\partial g(\tau)}{\partial \tau}$~\cite{SM}. Without
loss of generality, one may assume that the coordinate system is
defined on a square centred on the origin, as in
figure~\ref{fig:2}. For a point $p$ along the right-hand face of the
square above the $x$-axis, let $q$ be the corresponding point along
the orbit which intersects the top face, with the convention that when
$p=P$ lies on the stable manifold, the corresponding point is $q=Q$ on
the unstable manifold. The orbit consists of two segments: the segment
$\overline{pq}$ inside the box, and the segment $\overline{qp}$ lying
in the complement of the box. The period $T=T(p)$ of this orbit is the
sum of the time $T_0(p)$ that the orbit spends on the segment
$\overline{pq}$ plus the time $T_1(p)$ that the orbit spends on the
segment $\overline{qp}$. The time $T_1(p)$ is a real-analytic function
that approaches the finite limit $T_1(P)$ as $p \to P$; $T_0(p)$ is
also real-analytic and approaches $+\infty$ as $p \to P$.

From equation (\ref{eq:bnloc2}), one has the equation
\begin{equation} \label{eq:bnloc3}
\bn(p) = \frac{T_0}{T^2} \times \int_0^{T_0} \n \circ \psi^s_t(p)\, \d t
+ \frac{T_1}{T^2} \times \int_{T_0}^{T} \n \circ \psi^s_t(p)\, \d t.
\end{equation}
The second term is bounded by a constant times $\frac{T_1}{T}$, which
converges to $0$ as $p \to P$. The first term converges to
$\n(0)=\bn(0)=\bn(P)$ as $p \to P$.  

A similar, but slightly more involved, argument shows that if $p$ lies
in the right-hand face of the square below the $x$-axis, then $\bn(p)$
converges to $\bn(P)$, also. By symmetry and invariance of $\bn$ under
$\psi^s$, this proves that $\bn$ is a continuous function in a
neighbourhood of the hyperbolic singularity and its stable and
unstable manifold.

The reader may verify by direct computation that, if $\n = y$ in the
coordinate box, then $\frac{\partial \bn}{\partial y}$ diverges to
$+\infty$ as $p \to P$ ($y \to 0$).

\medskip
\noindent{\em Case 2, $|s|=1$:} In this case, $f|\Ss$ has two critical
points -- at $\a_0=\a_1=0,\n=\pm 1$ -- that are both degenerate. The
argument of case 1a may be adapted to show that $\bn$ is a continuous
function at each of these critical points.
\end{proof}

Here are some further properties of $\bn$. Since $\bn$ is
$\psi^s$-invariant, one may view it as a function defined on the image
of $f|\Ss$. In this case, it makes sense to say that $\bn$ is monotone
increasing.

\begin{Proposition}
If $s>1$ (resp. $s<-1$), then $\bn$ is a monotone increasing
(resp. decreasing) function that vanishes only on the zero level of
$f|\Ss$.
\label{PB}
\end{Proposition}

\begin{proof}
The symmetry of $f$ and the symplectic form $\theta$ dictate that
$\bn(c)$ be an odd function of $c$. Therefore $\bn$ always vanishes on
the zero level of $f$.

Let us suppose that $s>0$; the case where $s<0$ is analogous. From the
previous proposition, it suffices to prove that $\bn$ is monotone
increasing on the regular levels of $f|\Ss$. From equation
(\ref{eq:bnloc}), one sees that if $\n_2(\phi,I)>\n_1(\phi,I)$ for all
$\phi,I$, then $\bn_2(I) > \bn_1(I)$ for all $I$. For our purposes,
let $\n_1=\n$ and let $\n_2 = \n \circ \gamma_\tau$ where $\gamma$ is
a gradient-like flow for $f|\Ss$ -- that takes the form
$\gamma_\tau(\phi,I)=(\phi,I+\tau)$ in angle-action coordinates -- and
$\tau>0$ is a small positive number. That is, if the derivative of
$\n$ in the direction of the gradient-like flow $\gamma$ is positive,
then $\bn$ is a monotone increasing function. Let us remark that to
test the positivity of this directional derivative, it suffices to use
any gradient-like vector field; in particular, it suffices to compute
the directional derivative of $\n$ with respect to the standard
gradient vector field of $f|\Ss$.

One computes that 
\begin{equation} \label{eq:grad}
\left\langle \d \n , \nabla(f|\Ss) \right\rangle = \left[ 
\begin{array}{cc}
\a_0 & \a_1
\end{array}
\right]\,
\left[ 
\begin{array}{cc}
s & \n\\
\n & s
\end{array}
\right]\,
\left[ 
\begin{array}{c}
\a_0\\ \a_1
\end{array}
\right].
\end{equation}
The symmetric matrix is positive definite if $s>0$ and $s^2>\n^2$. If
$s>1$, then the matrix is always positive definite, whence the
right-hand side vanishes only at $\a_0=\a_1=0,\n=\pm 1$. This proves
the proposition.

\end{proof}

\begin{Remark}{\rm
When $|s|<1$, the function $\bar\n$ cannot be monotone increasing. As
one can see in figure \eqref{fig:1}, $\bar\n$ attains its maximum
value of unity at the hyperbolic fixed point $\a_0=\a_1=0,\n=1$; at
the same point $f=s$. On the other hand, at the elliptic critical
points $\a_0=\a_1=\pm\sqrt{\frac{1}{2}(1-s^2)}, \n=s$, $\bar\n$
attains a value of $s$ while $f=\frac{1}{2}(1+s^2)$. Thus: $s <
\frac{1}{2}(1+s^2)$ while $1=\bar\n(s) >
\bar\n(\frac{1}{2}(1+s^2))=s$. Numerical calculations do suggest that
$\bar\n$ is monotone increasing on $[-s,s]$ and decreasing on the two
complementary subintervals (see figure \ref{fig:3}).

A related issue concerns the monotone nature of the function $s
\mapsto h_{\leb}(\varphi^s)$. In figure \eqref{fig:4} we give evidence
from numerical computations that this function is a monotone function on
$(-\infty,0]$ and $[0,\infty)$. 

The function $\bn$ is approximated by integrating the Euler equations
in the almost canonical variables $\n,\phi$ (see the discussion around
equation \eqref{eq:theta}) using the Runge-Kutta $4$-step method and
averaging $\n$ over a numerically computed period. The function
$h_{\leb}(\varphi^s)$ is approximated by numerically integrating $\bn$
over a grid using Simpson's rule. Data and source code is available
from \href{http://www.maths.ed.ac.uk/~lbutler/dsol.html}{here}.

{
\def\Diamond{{\cdot}}
\def\Box{{\cdot}}
\begin{center}
\begin{figure}[htb]
\input{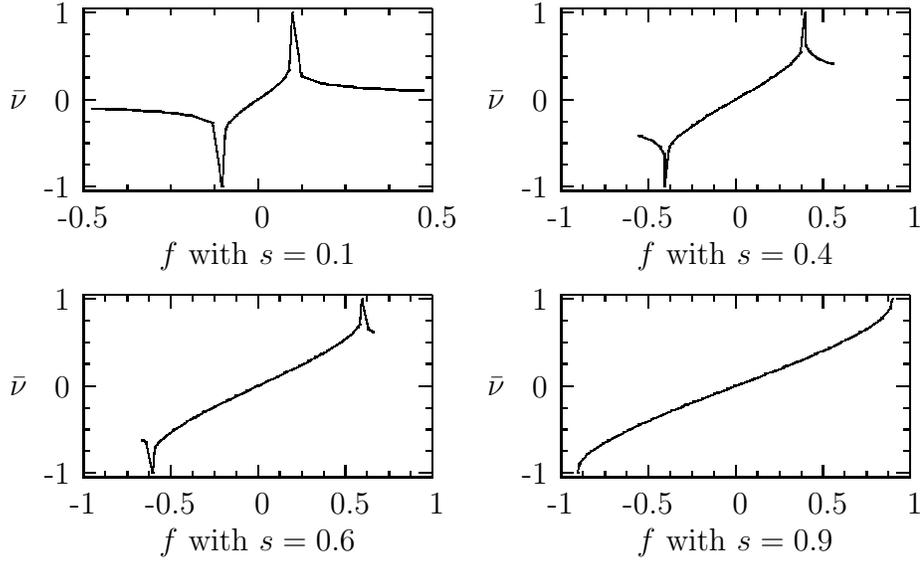}
\caption{The function $\bar\n$ as a function of $f$ for selected
  values of $s$. Note the loss of differentiability at the hyperbolic
  critical level $f=s$ and the lack of monotonicity. } \label{fig:3}
\end{figure}
\end{center}
}

{
\def\picA{\objectbox{
\setlength{\unitlength}{0.240900pt}
\ifx\plotpoint\undefined\newsavebox{\plotpoint}\fi
\sbox{\plotpoint}{\rule[-0.200pt]{0.400pt}{0.400pt}}%
\begin{picture}(1500,900)(0,0)
\sbox{\plotpoint}{\rule[-0.200pt]{0.400pt}{0.400pt}}%
\put(181.0,137.0){\rule[-0.200pt]{4.818pt}{0.400pt}}
\put(161,137){\makebox(0,0)[r]{ 0}}
\put(1419.0,137.0){\rule[-0.200pt]{4.818pt}{0.400pt}}
\put(181.0,279.0){\rule[-0.200pt]{4.818pt}{0.400pt}}
\put(161,279){\makebox(0,0)[r]{ 0.1}}
\put(1419.0,279.0){\rule[-0.200pt]{4.818pt}{0.400pt}}
\put(181.0,421.0){\rule[-0.200pt]{4.818pt}{0.400pt}}
\put(161,421){\makebox(0,0)[r]{ 0.2}}
\put(1419.0,421.0){\rule[-0.200pt]{4.818pt}{0.400pt}}
\put(181.0,562.0){\rule[-0.200pt]{4.818pt}{0.400pt}}
\put(161,562){\makebox(0,0)[r]{ 0.3}}
\put(1419.0,562.0){\rule[-0.200pt]{4.818pt}{0.400pt}}
\put(181.0,704.0){\rule[-0.200pt]{4.818pt}{0.400pt}}
\put(161,704){\makebox(0,0)[r]{ 0.4}}
\put(1419.0,704.0){\rule[-0.200pt]{4.818pt}{0.400pt}}
\put(181.0,846.0){\rule[-0.200pt]{4.818pt}{0.400pt}}
\put(161,846){\makebox(0,0)[r]{ 0.5}}
\put(1419.0,846.0){\rule[-0.200pt]{4.818pt}{0.400pt}}
\put(307.0,123.0){\rule[-0.200pt]{0.400pt}{4.818pt}}
\put(307,82){\makebox(0,0){-4}}
\put(307.0,840.0){\rule[-0.200pt]{0.400pt}{4.818pt}}
\put(558.0,123.0){\rule[-0.200pt]{0.400pt}{4.818pt}}
\put(558,82){\makebox(0,0){-2}}
\put(558.0,840.0){\rule[-0.200pt]{0.400pt}{4.818pt}}
\put(810.0,123.0){\rule[-0.200pt]{0.400pt}{4.818pt}}
\put(810,82){\makebox(0,0){ 0}}
\put(810.0,840.0){\rule[-0.200pt]{0.400pt}{4.818pt}}
\put(1062.0,123.0){\rule[-0.200pt]{0.400pt}{4.818pt}}
\put(1062,82){\makebox(0,0){ 2}}
\put(1062.0,840.0){\rule[-0.200pt]{0.400pt}{4.818pt}}
\put(1313.0,123.0){\rule[-0.200pt]{0.400pt}{4.818pt}}
\put(1313,82){\makebox(0,0){ 4}}
\put(1313.0,840.0){\rule[-0.200pt]{0.400pt}{4.818pt}}
\put(181.0,123.0){\rule[-0.200pt]{303.052pt}{0.400pt}}
\put(1439.0,123.0){\rule[-0.200pt]{0.400pt}{177.543pt}}
\put(181.0,860.0){\rule[-0.200pt]{303.052pt}{0.400pt}}
\put(181.0,123.0){\rule[-0.200pt]{0.400pt}{177.543pt}}
\put(40,491){\makebox(0,0){\ensuremath{h_{\leb}(\varphi^s)}}}
\put(810,21){\makebox(0,0){\ensuremath{s}}}
\put(181,846){\usebox{\plotpoint}}
\put(181.0,846.0){\rule[-0.200pt]{303.052pt}{0.400pt}}
\put(181,836){\usebox{\plotpoint}}
\multiput(181,836)(20.756,0.000){4}{\usebox{\plotpoint}}
\multiput(244,836)(20.756,0.000){3}{\usebox{\plotpoint}}
\multiput(307,836)(20.753,-0.329){3}{\usebox{\plotpoint}}
\multiput(370,835)(20.753,-0.329){3}{\usebox{\plotpoint}}
\multiput(433,834)(20.712,-1.336){3}{\usebox{\plotpoint}}
\multiput(495,830)(20.662,-1.968){3}{\usebox{\plotpoint}}
\put(573.45,822.00){\usebox{\plotpoint}}
\put(593.02,818.49){\usebox{\plotpoint}}
\put(612.00,814.00){\usebox{\plotpoint}}
\put(630.57,807.43){\usebox{\plotpoint}}
\put(649.14,800.86){\usebox{\plotpoint}}
\put(666.78,792.11){\usebox{\plotpoint}}
\put(683.49,781.00){\usebox{\plotpoint}}
\put(701.66,775.66){\usebox{\plotpoint}}
\put(715.03,761.46){\usebox{\plotpoint}}
\put(727.26,747.68){\usebox{\plotpoint}}
\put(736.32,729.68){\usebox{\plotpoint}}
\put(744.99,719.04){\usebox{\plotpoint}}
\put(751.87,706.38){\usebox{\plotpoint}}
\put(758.33,687.53){\usebox{\plotpoint}}
\put(762.55,667.49){\usebox{\plotpoint}}
\put(766.60,647.21){\usebox{\plotpoint}}
\put(771.94,628.00){\usebox{\plotpoint}}
\put(775.31,609.35){\usebox{\plotpoint}}
\put(776.86,588.66){\usebox{\plotpoint}}
\put(780.08,568.17){\usebox{\plotpoint}}
\put(781.88,547.49){\usebox{\plotpoint}}
\put(786.27,528.06){\usebox{\plotpoint}}
\put(788.67,507.45){\usebox{\plotpoint}}
\put(790.16,486.75){\usebox{\plotpoint}}
\put(791.69,466.05){\usebox{\plotpoint}}
\put(795.14,445.61){\usebox{\plotpoint}}
\put(797.77,425.03){\usebox{\plotpoint}}
\put(799.70,404.38){\usebox{\plotpoint}}
\put(801.26,383.70){\usebox{\plotpoint}}
\put(803.25,363.06){\usebox{\plotpoint}}
\put(804.48,342.35){\usebox{\plotpoint}}
\put(805.47,321.61){\usebox{\plotpoint}}
\put(806.53,300.89){\usebox{\plotpoint}}
\put(808.18,280.20){\usebox{\plotpoint}}
\put(809.38,259.49){\usebox{\plotpoint}}
\multiput(810,239)(0.629,20.746){2}{\usebox{\plotpoint}}
\put(811.87,280.71){\usebox{\plotpoint}}
\put(813.49,301.39){\usebox{\plotpoint}}
\put(814.56,322.12){\usebox{\plotpoint}}
\put(815.54,342.85){\usebox{\plotpoint}}
\put(816.78,363.57){\usebox{\plotpoint}}
\put(818.84,384.20){\usebox{\plotpoint}}
\put(820.33,404.89){\usebox{\plotpoint}}
\put(822.31,425.53){\usebox{\plotpoint}}
\put(824.91,446.11){\usebox{\plotpoint}}
\put(828.35,466.55){\usebox{\plotpoint}}
\put(829.88,487.25){\usebox{\plotpoint}}
\put(831.39,507.95){\usebox{\plotpoint}}
\put(833.78,528.56){\usebox{\plotpoint}}
\put(838.17,548.00){\usebox{\plotpoint}}
\put(839.97,568.67){\usebox{\plotpoint}}
\put(843.18,589.16){\usebox{\plotpoint}}
\put(844.72,609.86){\usebox{\plotpoint}}
\put(848.57,628.00){\usebox{\plotpoint}}
\put(853.46,647.72){\usebox{\plotpoint}}
\put(857.50,667.99){\usebox{\plotpoint}}
\put(861.78,688.03){\usebox{\plotpoint}}
\put(868.17,706.89){\usebox{\plotpoint}}
\put(875.17,719.52){\usebox{\plotpoint}}
\put(884.04,730.04){\usebox{\plotpoint}}
\put(892.84,748.18){\usebox{\plotpoint}}
\put(905.25,761.88){\usebox{\plotpoint}}
\put(918.69,775.31){\usebox{\plotpoint}}
\put(937.01,781.01){\usebox{\plotpoint}}
\put(953.67,792.34){\usebox{\plotpoint}}
\put(971.31,801.00){\usebox{\plotpoint}}
\put(989.79,807.79){\usebox{\plotpoint}}
\put(1008.51,814.00){\usebox{\plotpoint}}
\put(1027.44,818.72){\usebox{\plotpoint}}
\put(1047.06,822.00){\usebox{\plotpoint}}
\multiput(1062,824)(20.659,1.999){3}{\usebox{\plotpoint}}
\multiput(1124,830)(20.714,1.315){3}{\usebox{\plotpoint}}
\multiput(1187,834)(20.753,0.329){3}{\usebox{\plotpoint}}
\multiput(1250,835)(20.753,0.329){3}{\usebox{\plotpoint}}
\multiput(1313,836)(20.756,0.000){3}{\usebox{\plotpoint}}
\multiput(1376,836)(20.756,0.000){3}{\usebox{\plotpoint}}
\put(1439,836){\usebox{\plotpoint}}
\put(181.0,123.0){\rule[-0.200pt]{303.052pt}{0.400pt}}
\put(1439.0,123.0){\rule[-0.200pt]{0.400pt}{177.543pt}}
\put(181.0,860.0){\rule[-0.200pt]{303.052pt}{0.400pt}}
\put(181.0,123.0){\rule[-0.200pt]{0.400pt}{177.543pt}}
\end{picture}
}}
\def\picB{\begin{rotate}{-90}\resizebox{3cm}{!}{\includegraphics[width=4cm, height=4cm]{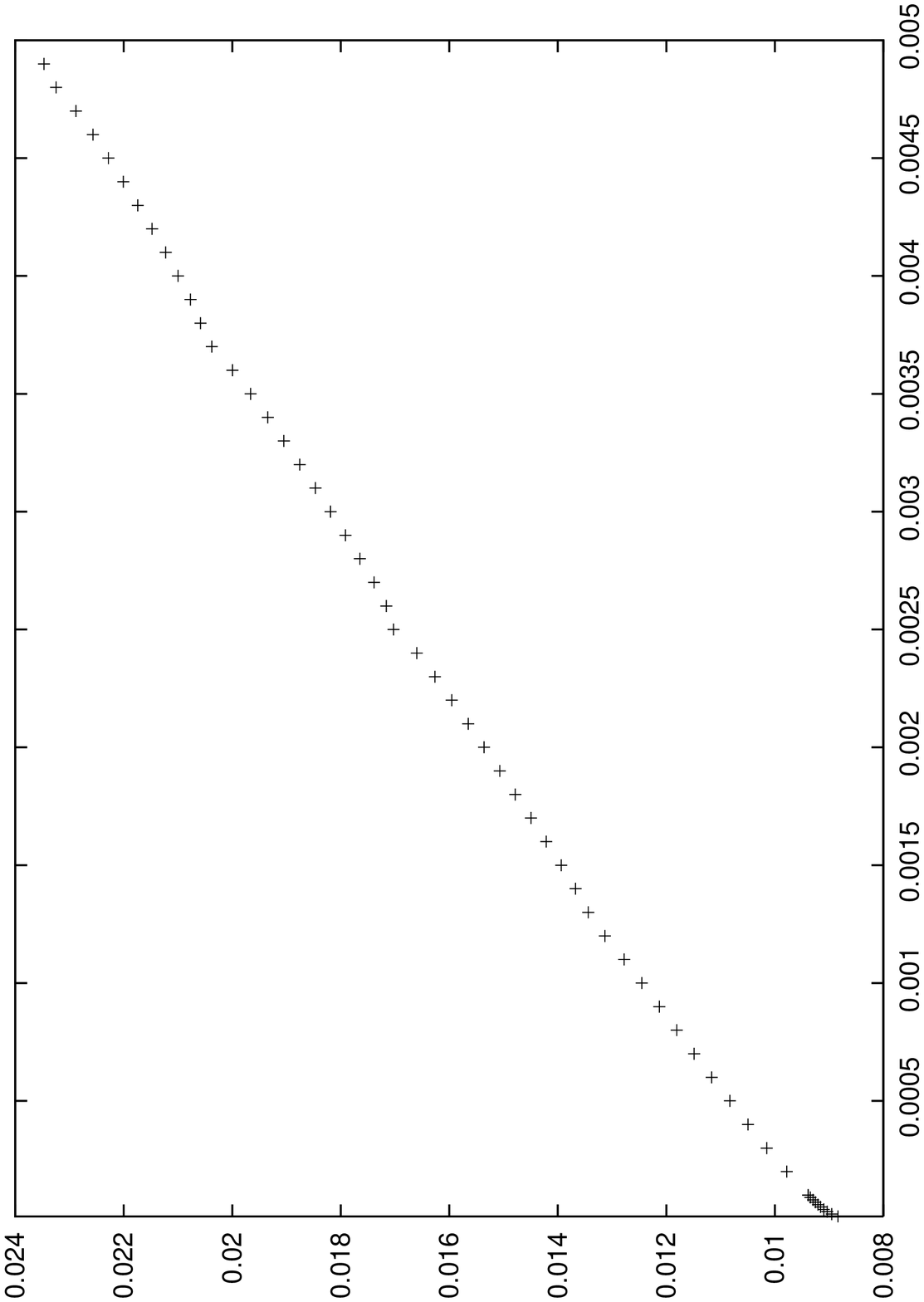}}\end{rotate}}
\def\picC{\begin{rotate}{-90}\resizebox{3cm}{!}{\includegraphics[width=4cm, height=4cm]{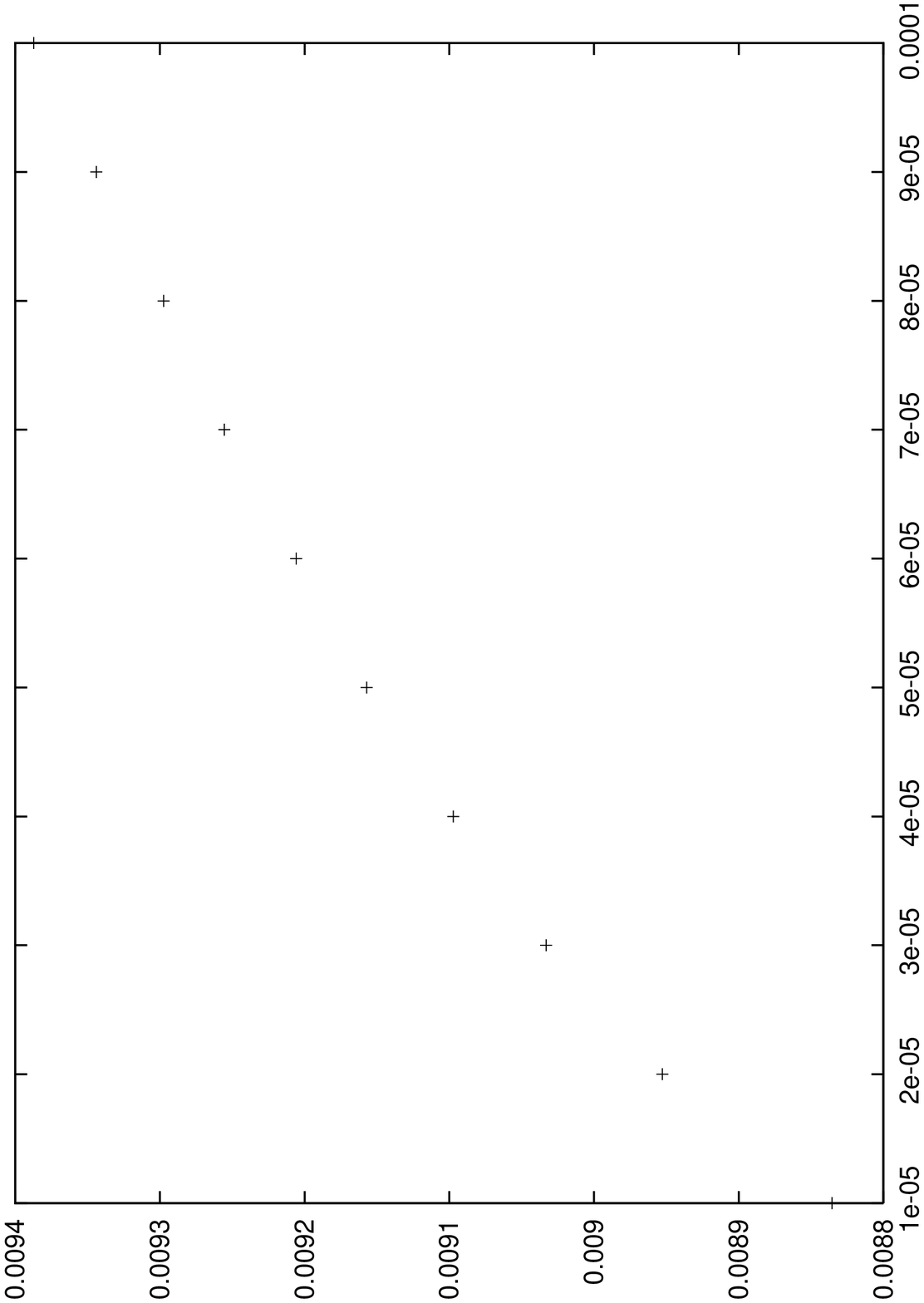}}\end{rotate}}

\begin{center}
\begin{figure}[htb]
\begin{xy}
(0,0)*{\objectbox{\input{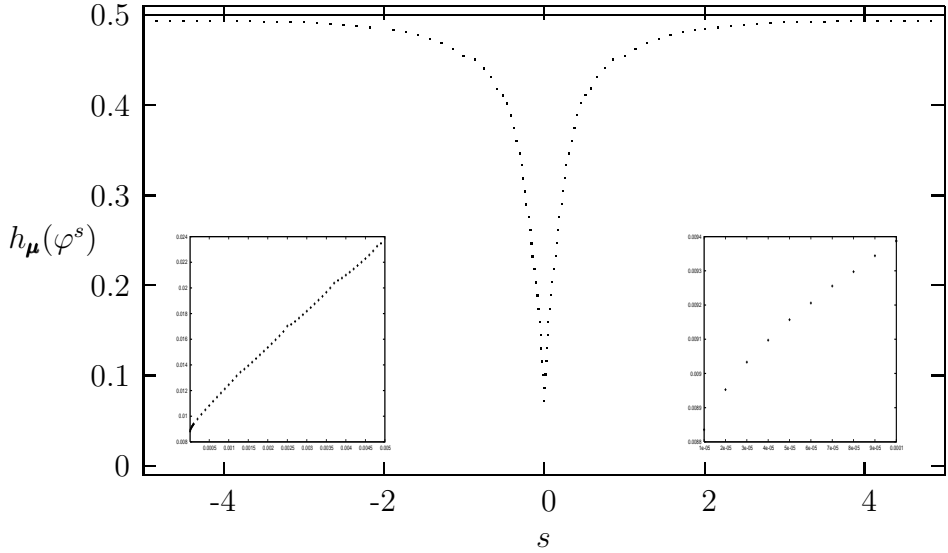}}} ;
(23,5)*{\picC} ;
(-45,5)*{\picB} ;
\end{xy}
\caption{The function $h_{\leb}(\varphi^s)$ as a function of $s$. Inset (left):
on the interval $[0, 5 \times 10^{-3}]$; Inset (right): on the
interval $[0, 1 \times 10^{-4}]$.  } \label{fig:4}
\end{figure}
\end{center}
}
}
\end{Remark}

\begin{Remark}{\rm Consider an orbit of the magnetic flow on $\Sol$ that
projects onto a closed orbit of $E_h$.  From equation (\ref{eq:XH}) it
is clear that $u$ is a periodic function of time if and only if
$\bar{\n}=0$.  Left-invariance--or an easy check using
(\ref{eq:XH})--gives that the functions $p_{y_{0}}+sy_1$ and
$p_{y_{1}}-sy_0$ are first integrals in $\Sol$.  Since $\alpha_0$ and
$\alpha_1$ are periodic, we conclude that $p_{y_{0}}=e^{-u}\alpha_{0}$
and $p_{y_{1}}=e^{u}\alpha_{1}$ are periodic if $u$ is periodic. Thus,
if $s>0$ and $\bar{\n}=0$, the orbit of the magnetic flow on $\Sol$
is periodic. Since there are always closed orbits of $E_h$ with
$\bar{\n}=0$ we conclude that for $s>0$ the magnetic flow on
$\DSol\setminus \S$ always has contractible closed orbits.

Observe that for the geodesic flow ($s=0$) no closed orbit is contractible,
since if $u$ is periodic, $y_0$ and $y_1$ must diverge linearly.
}
\label{clor}
\end{Remark}

\subsection{Cocompact subgroups of $\Sol$} To compute the metric
entropy of the magnetic flow, it is useful to view the lattice
subgroup $\DSol$ of $\Sol$, especially the diagonalizing
transformation $P$ described in the introduction, intrinsically.

Given a lattice subgroup $\DSol$ of $\Sol$, there is an exact sequence
$\Z^2 \hookrightarrow \DSol \to \Z$ induced by the exact sequence
$\R^2 \hookrightarrow \Sol \to \R$ \cite[pp. 470--472]{Scott}. The
quotient group $\Z$ acts on $\Z^2$ via a representation $\rho : \Z \to
SL(2,\Z)$. The generator $\rho(1)$ ($=A$ from the introduction) is a
hyperbolic matrix with eigenvalues $\lambda^{\pm 1}$, $| \lambda |>1$.
In terms of the coordinate system (equation \ref{eq:sm}), the group
$\DSol$ can be described as follows. Let $F=\Q(\lambda)$ be the
quadratic number field obtained by adjoining $\lambda$ to the
rationals. The integers of $F$, $\O$, is isomorphic to $\Z^2$ as an
abelian group, and the unit group of $\O$, $\U$, acts as an
automorphism group. The group $\DSol$ is naturally isomorphic to a
finite-index subgroup of the semi-direct product $\U \star \O$. We
shall henceforth identify $\DSol$ with a subgroup of $\U \star \O$.

The volume of $\DSol \lhd \U \star \O$ can be defined to be 
\begin{equation} \label{eq:volume}
\vol{\DSol} := \log|\lambda| \times \det \left[ \begin{array}{cc}
a_0^{(0)} & a_1^{(0)}\\
a_0^{(1)} & a_1^{(1)}
\end{array}\right],
\end{equation}
where $a_0,a_1$ generate $\DSol \cap \O$ and $a_i^{(j)}$ is the $j$-th
conjugate of $a_i$%
\footnote{The field $\Q(\lambda)$ is a quadratic number field and so
  equals $\Q(\sqrt{d})$ for some positive, square-free integer
  $d$. The map $\sqrt{d} \mapsto -\sqrt{d}$ induces a field
  automorphism, and the image of $a=a^{(0)}$ under this automorphism
  is refered to as a conjugate of $a$ and denoted by $a^{(1)}$.}. One
  can see that $\vol{\DSol}$ is the determinant of the injection of
  $\Z^2 \ltimes_A \Z$ into $\Sol$ defined in the introduction; indeed,
  the matrix $P$ introduced there is effectively the matrix on the
  right-hand side of equation (\ref{eq:volume}). It is clear that
  $\vol{\DSol}$ is the volume of a fundamental region for $\DSol$ in
  $\Sol$ relative to the volume form $\d\u \wedge \d\y_0 \wedge
  \d\y_1$. That is
\begin{equation} \label{eq:vol}
\vol{(\DSol \backslash \Sol)} = \vol{\DSol}.
\end{equation}
Let $\Sigma = \DSol\backslash \Sol$ and let $\leb$ be the $X_H$-invariant
probability measure on $\Sigma \times \Ss$ induced by $\omega_s^3 = -
\d\u \wedge \d \y_0 \wedge \d \y_1 \wedge \d \n \wedge \d \a_0 \wedge
\d \a_1$, {\em i.e.}
$$\leb = \frac{1}{\vol{\DSol}} \times \d\u \wedge \d \y_0 \wedge \d
\y_1 \wedge \theta.$$

\subsection{Metric Entropy of the Magnetic Flow}
$\;$

\medskip

\noindent{\bf Theorem A.}
{\it Let $s\neq0$ and $\varphi^s : \R \times \Sigma \times \Ss \to \Sigma \times \Ss$ be
the magnetic flow with infinitesimal generator $X_H$. The
metric entropy of the time-$1$ map $\varphi^s_1$ is
\begin{equation} \label{eq:entropy}
h_{\leb}(\varphi^s_1) = \int_{\Ss}\ |\bn|\, \d\theta.
\end{equation}
Therefore, since $\bn$ is non-zero on a positive measure set,
$h_{\leb}(\varphi^s_1)>0$.  Moreover, $h_{\mu}(\varphi^{s}_1)$
approaches $1/2$ as $s\to\infty$ and $h_{top}(\varphi^s_1)\equiv 1$.
}

\medskip

\begin{Remark} {\rm As is proven in Proposition \ref{PA},
 $\bn$ is a continuous function that is real-analytic on the
    complement of the critical levels of $f|\Ss$. Therefore $\bn$ is
    non-zero on a set of full measure. It is almost certain that $\bn$
    vanishes only on one level of $f|\Ss$; Proposition \ref{PB} proves this
    when $|s|>1$.}
\end{Remark}

\begin{proof}
First, consider a flow $\varphi : \R \times \Sol \times \T^1 \to \Sol
\times \T^1$ which is a skew product over a translation
\begin{equation} \label{eq:skew}
\varphi_t(g,\phi) = (\cc(t,g,\phi),\phi+at \bmod 1) \hspace{20mm} \forall g\in\Sol, \phi\in\T^1.
\end{equation}
If $\varphi$ is assumed to be left-invariant, then the cocycle $\cc$
satisfies $\cc(t,g,\phi) = g\star \cc(t,1,\phi)$. Therefore, if $T = 1/a$, then
\begin{equation} \label{eq:skew2}
\varphi_T(g,\phi) = (g \star \cc(\phi), \phi \bmod 1) \hspace{20mm} \forall g\in\Sol, \phi\in\T^1,
\end{equation}
where $\cc(\phi)=\cc(T,1,\phi)$. Therefore, for all $n\in\Z$,
\begin{equation} \label{eq:skew3}
\varphi_{nT}(\DSol g,\phi) = (\DSol g \star \cc(\phi)^n, \phi \bmod 1) \hspace{20mm}
\forall g\in\Sol, \phi\in\T^1,
\end{equation}
The cocycle $\cc(\phi) \in \Sol$ either takes values in the non-hyperbolic
subgroup $\R^2$ or it has a non-trivial projection to $\R$. In the
former case, $\varphi_{T}$ has zero entropy. In the latter case,
$T\Sol$ splits into $3$ complementary, left-invariant line-bundles
$E^+, E^-$ and $E^0$. These line bundles are determined by their value
at the identity of $\Sol$. If we identify $T_I \Sol$ as the Lie
algebra $\sol$, then $E^+$ is the unstable subspace, $E^-$ is the
stable subspace and $E^0$ is the centralizer of ${\rm Ad}_{\cc(\phi)}$,
respectively. Let $\lambda_+(g)$ be the $\log$ of the largest
eigenvalue of ${\rm Ad}_g$, $g\in\Sol$. One sees using Pesin's formula that 
\begin{equation} \label{eq:hmu}
h_{\leb_c}(\varphi_1) = \frac{1}{T} \times \int_0^1
\lambda_+(\cc(\phi))\, \d \phi,
\end{equation}
where $\leb_c = \frac{1}{\vol{\DSol}} \times \d\u \wedge \d \y_0
\wedge \d \y_1 \wedge \d \phi$ is a $\varphi$-invariant probability
measure on $\Sigma \times \T^1$.

A simple computation shows that $\lambda_+(g)$ is the projection $g
\mapsto |u(g)|$ induced by $\Sol \stackrel{u}{\longrightarrow} \R$. In
addition, if we observe that $\cc(\phi) = \cc(0,1,\phi)^{-1} \star
\cc(T,1,\phi)$ and use the fact that $u$ is a group homomorphism, then
\begin{equation} \label{eq:du0}
|\Delta u| = |u( \cc(\phi) )| = \lambda_+( \cc(\phi) ),
\end{equation}
where $\Delta u$ is the change in $u$ over the time interval $[0,T]$.

\medskip
Let us turn to the magnetic flow: Let $c$ be a regular value of
$f|\Ss$ and introduce action-angle variables in a neighbourhood of
$\f_c \subset \Ss$. The flow, $\varphi^s$, of $X_H$ restricted to
$\Sigma \times \f_c$ is of the form described by equation
(\ref{eq:skew}), with $a = \frac{\partial g}{\partial I}$, see
equation (\ref{eq:m}). The Liouville measure $\leb$ on $\Sigma \times
\Ss$ induces the invariant conditional probability measure $\leb_c$ on
$\Sigma \times \f_c$. Inspection of equation (\ref{eq:XH}) shows that
over the period $T$, $\u$ changes by
\begin{equation} \label{eq:du}
\Delta u = \int_0^T\ \n(t)\, \d t.
\end{equation}
Since $\n$ is a periodic function, the integral for $\Delta u$ is
independent of the angle variable $\phi$. Equations
(\ref{eq:hmu}--\ref{eq:du0}) therefore show that $|\Delta u|/T$ is the
metric entropy of $\varphi^s | \Sigma \times \f_c$ with respect to the
conditional probability measure $\leb_c$. Using equation (\ref{eq:EhI}) in
action-angle coordinates, one obtains
\begin{equation} \label{eq:duI}
\frac{\Delta u}{T} = \frac{1}{T} \times \int_0^1\ \n(\phi,I)\, \d \phi \times
  \frac{\partial I}{\partial g} = \int_0^1
  \n(\phi,I)\, \d\phi = \bn,
\end{equation}
since $T=\frac{\partial I}{\partial g}$. Therefore, we can integrate
to obtain the metric entropy of the magnetic flow on the unit-sphere
bundle
\begin{equation} \label{eq:hmufc}
h_{\leb}(\varphi^s) = \int_{\Sigma \times \Ss} |\Delta u|\, \d\leb =
\int_{\Ss} |\bn|\, \d\theta.
\end{equation}
This proves equation (\ref{eq:entropy}). Let us prove the remaining
two points in Theorem A.

\subsubsection*{Topological entropy}
We note that the arguments above also imply that the sum of the
non-negative Liapunov exponents of $\varphi^s$ is given by
$|\bar{\n}|\leq 1$.  Thus by Ruelle's inequality and the variational
principle for topological entropy we see that $h_{top}(\varphi^s)\leq
1$. Since the flow on the set $p_u=\pm 1$, $p_{y_{0}}=p_{y_{1}}=0$ is
the same for all $s$ and carries entropy $1$ we conclude that
$h_{top}(\varphi^s)\equiv 1$.

\subsubsection*{The limit of metric entropy}
To compute $\lim_{s \to \infty} h_{\leb}(\varphi^s)$, note that
$\frac{1}{s} \times f = \n + \a_0\a_1 \times \frac{1}{s}$. Therefore,
as $s \to \infty$, the regular level sets of $f|\Ss$ converge
uniformly in the $C^1$ topology to the level sets of $\n$, {\it i.e.}
the regular level sets of $f|\Ss$ converge to circles at a constant
height off the $\a_0-\a_1$ plane (and $f|\Ss$ has only the points
$\a_0=\a_1=0,\n=\pm 1$ as critical points). 

Let us coordinatize $\Ss-\left\{ (0,0,\pm 1) \right\}$
by spherical coordinates
$$\a_0=\cos(2\pi\xi) \sin(\eta), \a_1=\sin(2\pi\xi)\sin(\eta),
\n=\cos(\eta) \qquad 0 < \eta < \pi, 0 \leq \xi < 1.$$ The angle
$\xi$ is the normalized longitudinal angle which vanishes along
$\left\{ \a_1=0,\a_0>0 \right\}$ and has $\frac{\partial \xi}{\partial
\a_1} > 0$ along the same privileged longitude. In spherical
coordinates, the normalized area form is
\begin{equation} \label{eq:theta}
\theta = \frac{1}{2} \times \sin(\eta)\, \d \eta \wedge \d \xi.
\end{equation}

For $s>1$, we normalize the action-angle coordinates
$(I,\phi)=(I_s,\phi_s)$ on $\Ss-\left\{ (0,0,\pm 1) \right\}$ as
follows: first, $\phi_s = 0$ and $\frac{\partial \phi_s}{\partial
\a_1} > 0$ along the privileged longitude $\left\{ \a_1=0,\a_0>0
\right\}$; second, $I_s(\mu)$ is defined to be the area of the
sublevel set $\left\{ f \leq f(\mu) \right\}$ in $\Ss$.  The above
paragraph shows that as $s \to \infty$, $I_s$ converges to the
function $I_{\infty}$ which gives the area of the region in $\Ss$
below height $\n$. A computation shows that
$I_{\infty}=\frac{1+\n}{2}$. On the other hand, $\phi_s$ converges to
the normalized longitudinal angle $\xi$. 

Inspection of equations (\ref{eq:m},\ref{eq:bnloc}) shows that the
mean value of $\n$ averaged with respect to the measure $\d \phi_s$
converges to $\n$ as $s \to \infty$. This convergence is in the
uniform $C^0$ topology. Therefore
\begin{equation} \label{eq:hlebs}
\lim_{s \to \infty} h_{\leb}(\varphi^s_1) = \int_0^1 \int_0^\pi
|\cos(\eta)\sin(\eta)|\, \d \eta\, \d \xi \times \frac{1}{2} = \frac{1}{2},
\end{equation}
as asserted.
\end{proof}

\subsection{A variation} There is an interesting variation of the
previous example. Consider the group $G=\S\times\re$ and the left-invariant
2-form given by $\Omega:=du\wedge dt$, where $t$ denotes the variable on the
$\re$-factor. We consider on $G$ the left-invariant metric given by
$\ds^2+dt^2$ and the cocompact lattice $\DSol\times\Z$.
The magnetic flow $\varphi^s$ on the compact quotient thus obtained
has the following remarkable
properties for $s\neq 0$ (as before $s$ is the intensity):
\begin{itemize}
\item $h_{top}(\varphi^s)=0$ for $s\neq 0$. This shows that topological
entropy may be discontinuous when a twist in the symplectic structure
is introduced;
\item $\varphi^s$ is completely integrable with {\it real analytic} integrals.
If we let $\tau:=p_t$, then the integrals are
 $\alpha_0\,\alpha_1$, $\alpha_0 e^{-\tau/s}$ (the two Casimirs) and
$\tau-su$, which can be made invariant under
the lattice just by composing with a suitable periodic function.
\end{itemize}
We leave the details of the proofs of these claims to the reader, but they do
follow in a straightforward fashion from an analysis quite similar to the one
done in this section.

\section{Proof of Theorem B}
\label{tb}

We first prove the following easy lemma.

\begin{Lemma} Let $\mathfrak g$ be a Lie algebra such that
$\mbox{\rm Ker}\,L$ is generated
by 2-vectors. Let $\Omega$ be an antisymmetric bilinear form on $\mathfrak g$
such that $\Omega(x,y)=0$ for all $x,y$ with $[x,y]=0$.
Then $\Omega$ is exact, that is, there exists $b\in \mathfrak g^{*}$ such that
$\Omega(x,y)=b([x,y])$ for all $x,y\in\mathfrak g$.
\label{easy}
\end{Lemma}

\begin{proof} Let $L^*:\mathfrak g^*\to (\Lambda_{2}(\mathfrak g))^*$ be the dual of $L:\Lambda_{2}(\mathfrak g)\to\mathfrak g$. It suffices to show that
$\Omega$ is in the image of $L^*$. But the image of $L^*$ coincides with
the annihilator of $\mbox{\rm Ker}\,L$, so it suffices to check that
$\Omega(q)=0$ for all $q\in \mbox{\rm Ker}\,L$. But if $\mbox{\rm Ker}\,L$
is generated by 2-vectors we may write $q=\sum_{i}x_{i}\wedge y_{i}$ with
$x_{i}\wedge y_{i}\in \mbox{\rm Ker}\,L$. Thus
$\Omega(q)=\sum_i\Omega(x_i,y_i)$. But since $[x_i,y_i]=0$, $\Omega(x_i,y_i)=0$
by hypothesis and $\Omega(q)=0$.

\end{proof}

We now break the proof of Theorem B into a few simple steps.

\begin{enumerate}

\item Let $\sigma$ be a closed 2-form in $\Sigma$ with non-zero
  cohomology class. By a theorem of A. Hattori \cite{H} (which in turn
  is a generalization of a theorem of K. Nomizu for nilmanifolds
  \cite{N}), there exists a left-invariant closed 2-form $\Omega$
  cohomologous to $\sigma$. This is the only part of the proof in
  which we use that $G$ is completely solvable.  We denote by the same
  symbol $\Omega$ the 2-form on $G$ or on $\Sigma$.

 Write $\sigma=\Omega+d\theta$ for some smooth 1-form $\theta$.  The
 fibrewise shift $(x,p)\mapsto (x,p-\theta)$ takes compact sets to
 compact sets and is a symplectomorphism between
 $(T^*\Sigma,\omega_{\sigma})$ and $(T^*\Sigma,
 \omega_{\Omega})$. Hence, from now on we may suppose that the
 monopole is given by a closed left-invariant 2-form $\Omega$.

\item We identify $T^*G$ with $G\times\g^*$ using left translations.
Smooth left-invariant functions on $T^*G$ are then identified
with $C^{\infty}(\g^*)$. The twisted symplectic structure $\omega_{\Omega}$
determines a Poisson bracket $\{\;,\;\}_{\Omega}$.
Given $f,g\in C^{\infty}(\g^*)$ we have
\begin{equation}
\{f,g\}_{\Omega}(m)=m([d_{m}f,d_{m}g])+\Omega(d_{m}f,d_{m}g)
\label{fp}
\end{equation}
for every $m\in\g^*$ where $d_{m}f,d_{m}g\in\g$ using the canonical
isomorphism $(\g^{*})^{*}=\g$. This formula is a simple
consequence of the definition of the twisted symplectic
form on $T^*G$ plus left-invariance.

\item If $f,g\in C^{\infty}(\g^*)$ then, they induce functions
on $T^{*}\Sigma=\Sigma\times \g^*$, which only depend
on the $\g^*$-variables and their Poisson brackets is computed, of course, also using
(\ref{fp}).

\item Since the cohomology class of $\Omega$ is not zero, we now invoke
Lemma \ref{easy} to obtain two vectors $x,y\in\g$ such that
$[x,y]=0$ but $\Omega(x,y)\neq 0$.

\item The vectors $x,y$ are obviously linearly independent. Consider
a basis $\{e_1=x,e_2=y,e_3,\dots,e_n\}$ of $\g$ and let
$\{e_{1}^*,\dots,e_{n}^*\}$ be its dual basis. Given $m\in\g^*$
write $m=\sum_i m_i\,e_{i}^*$ and let $f_{i}(m):=m_i$.
Using (\ref{fp}) we have
\[\{f_{1},f_{2}\}_{\Omega}(m)=\Omega(x,y)\neq 0.\]

\item Consider $f_{1}$ as Hamiltonian on $T^{*}\Sigma$.
Along the Hamiltonian flow of $f_{1}$ we have
\[\dot{m}_{2}=\Omega(x,y)\neq 0\]
which readily implies that any compact set in $T^{*}\Sigma$ may be
displaced using the Hamiltonian flow of a suitable cut-off of $f_1$.
This finishes the proof of Theorem B.

\end{enumerate}

\subsection{Examples} Consider the Heisenberg Lie algebra $\mathfrak h_{2n+1}$
with basis 
$$\{x_{1},\dots,x_{n},y_{1},\dots,y_{n},z\}$$
 and non-zero brackets $[x_{i},y_{i}]=z$ for $i=1,\dots,n$. The image of
$L:\Lambda_{2}(\mathfrak h_{2n+1})\to \mathfrak h_{2n+1}$ is obviously one dimensional
and generated by $z$. All the vectors $x_i\wedge x_j$, $y_i\wedge y_{j}$, $x_i\wedge z$,
$y_i\wedge z$ and $x_i\wedge y_{j}$ for $i\neq j$ are in the Kernel of $L$.
Additional $n-1$, 2-vectors in the Kernel of $L$ are given by
\[(x_{i}+y_{1})\wedge(x_{1}+ y_{i})\]
for $i=2,\dots,n$. Thus $\mbox{\rm Ker}\,L$ is generated
by 2-vectors.

Another well known nilpotent Lie algebra $\mathfrak g_{2n+1}$ is given by a basis 
$$\{x_{1},\dots,x_{n},y_{1},\dots,y_{n},z\}$$
and non-zero brackets $[z,x_{i}]=y_i$ for $i=1,\dots,n$.
Here the kernel of $L$ is generated by $x_{i}\wedge x_j$, $y_{i}\wedge y_{j}$, $x_{i}\wedge y_{j}$
and $z\wedge y_{i}$.

Finally consider the nilpotent Lie algebra $\mathfrak u_{n}$ of upper
triangular $n\times n$ matrices (with zeros along the diagonal).  If
$e_{ij}$ denotes the matrix which has a $1$ in its $(i,j)$-entry and
zero everywhere else, then the non-zero brackets are
$[e_{ij},e_{jl}]=e_{il}$ where $i<j<l$. As in the case of the
Heisenberg Lie algebra it is easy to check that $\mbox{\rm Ker}\,L$ is
generated by 2-vectors and we leave this to the reader.

The simply connected nilpotent Lie groups associated with these Lie
algebras admit cocompact lattices and monopoles. The corresponding
second Betti numbers are:
\begin{align*}
b_{2}(\mathfrak h_{2n+1})&=2n^2-n-1,\;\;\;n\geq 2;\\
b_{2}(\mathfrak g_{2n+1})&=n(n+1);\\
b_{2}(\mathfrak u_{n})&=\frac{(n-2)(n+1)}{2}.\\
\end{align*}

To all of them Theorem B applies.

\subsection{Relation with Ma\~n\'e's critical value}
\label{mcv}

Let $M$ be a closed manifold and $\sigma$ a non-zero closed 2-form. We say
that a compact set $K\subset (T^{*}M,\omega_{\sigma})$ is
{\it stably displaceable} if $K\times S^1$ is displaceable
in $(T^{*}M\times T^{*}S^1,\omega_{\sigma}\oplus\omega_{0})$.
Let $g$ be a Riemannian metric on $M$.
Following Schlenk in \cite{S} we define
$d(g,\sigma)$ as the supremum of the values of $k\in\R$ such that
the set of $(x,p)\in T^*M$ with $|p|_{x}^2\leq 2k$ is stably displaceable.

The results of Laudenbach-Sikorav \cite{LS} and Polterovich \cite{Po}
that we mentioned in the Introduction imply that $d(g,\sigma)>0$. We
have introduced stable displacement to include the case in which the
Euler characteristic of $M$ is different from zero. Note that this was
unnecessary before because all the manifolds we discussed had vanishing
Euler characteristic.

Suppose now that $\sigma$ is {\it weakly exact}, that is, its lift
$\widetilde{\sigma}$ to the universal covering $\widetilde{M}$ of $M$
is exact.  Ma\~n\'e's critical value $c(g,\sigma)$ is defined as
\cite{BP}:
\[c(g,\sigma):=\inf_{u\in C^{\infty}(\widetilde{M},\re)}\;\sup_{x\in \widetilde{M}}\;
   \frac{1}{2}|d_{x}u+\theta_{x}|^{2},\]
where $\theta$ is any primitive of $\widetilde{\sigma}$.
As $u$ ranges over $C^{\infty}(\widetilde{M},\re)$ the form $\theta+du$ 
ranges over all primitives of $\widetilde{\sigma}$, because any two
primitives differ by a closed 1-form which must be exact since $\widetilde{M}$
is simply connected.
The critical value $c(g,\sigma)<\infty$ if and only if $\widetilde{\sigma}$
has bounded primitives.

\medskip

\noindent{\bf Question.} Is $d(g,\sigma)=c(g,\sigma)$ always?

\medskip

As far as we are aware, there are no counterexamples to this equality
which is motivated by the desire to relate Aubry-Mather theory
with Symplectic Topology. A full motivation for this question together
with more examples where equality holds maybe found in \cite{CFP}.
Suppose that $\pi_{1}(M)$ is amenable. Then (see \cite[Corollary 5.4]{P})
$c(g,\sigma)=\infty$ if and only if $[\sigma]\neq 0$.
Thus, to test the Question when $\pi_1(M)$ is amenable
and $\sigma$ is a monopole, we must show that $d(g,\sigma)=\infty$.
This is exactly the content of Theorem B which could then be interpreted
as evidence of a positive answer to the Question (recall that solvable groups are amenable).

\end{document}